\documentclass[12pt]{amsart}
\usepackage{graphicx, amssymb}
\numberwithin{equation}{section}

\newtheorem{Thm}{Theorem}[section]

\newtheorem{Lem}{Lemma}[section]

\def\N {\mathbf{N}}
\def\R {\mathbf{R}}
\def\S {\mathbf{S}}

\def\DD {\mathcal{D}}
\def\MM {\mathcal{M}}

\def\a {{\alpha}}
\def\de {{\delta}}
\def\th {{\theta}}
\def\l {{\lambda}}
\def\si {{\sigma}}
\def\om {{\omega}}

\def\rstr {{\big |}}
\def\indc {{\bf 1}}

\def\la {\langle}
\def\ra {\rangle}

\def\d {{\partial}}
\newcommand{\Div}{\operatorname{div}}
\def\grad {{\nabla}}
\newcommand{\curl}{\operatorname{curl}}
\newcommand{\Supp}{\operatorname{supp}}
\newcommand{\VP}{\operatorname{vp}}

\begin{document}
\title[Classical Solutions to Vlasov-Maxwell]
      {On Classical Solutions to the 3D\\
        Relativistic Vlasov-Maxwell System:\\
        Glassey-Strauss' Theorem Revisited}

\date{\today}

\author[F. Bouchut]{Fran\c cois Bouchut}
\address[F. B.]%
{CNRS, \& D\'epartement de Math\'ematiques et Applications\\
 Ecole Normale Sup\'erieure\\
 45 rue d'Ulm, F75230 Paris cedex 05}
\email{fbouchut@dma.ens.fr}

\author[F. Golse]{Fran\c cois Golse}
\address[F. G.]%
{Institut Universitaire de France\\
\& Laboratoire Jacques-Louis Lions\\
 Universit\'e Paris VI\\
 175 rue du Chevaleret, F75013 Paris}
\email{golse@math.jussieu.fr}

\author[C. Pallard]{Christophe Pallard}
\address[C. P.]%
{D\'epartement de Math\'ematiques et Applications\\
 Ecole Normale Sup\'erieure\\
 45 rue d'Ulm, F75230 Paris cedex 05}
\email{pallard@dma.ens.fr}

\begin{abstract}
R. Glassey and W. Strauss have proved in [Arch. Rational Mech. Anal. 92 (1986), 59--90]
that $C^1$ solutions to the relativistic Vlasov-Maxwell system in three space dimensions
do not develop singularities as long as the support of the distribution function in the
momentum variable remains bounded. The present paper simplifies their proof.
\end{abstract}

\maketitle

\section{Introduction}\label{INTRO}

\subsection{The relativistic Vlasov-Maxwell system}

The Vlasov-Max- well system is a mean-field, kinetic model for plasmas. Within the formalism
of kinetic theory, it describes the motion of a gas of charged, relativistic particles (e.g.
electrons or ions). Each particle is subject to the electromagnetic field created by all the
other particles, but not to its own self-consistent electromagnetic field which is neglected
in this model.

For simplicity, we give up the constraint of global neutrality and consider only the case of
a single species of charged particles, with distribution function denoted by $f$. Precisely,
$f(t,x,\xi)$ is the phase space density of particles which at time $t>0$ are located at the
point $x\in\R^3$ and have momentum $\xi\in\R^3$. Let $E\equiv E(t,x)$ and $B\equiv B(t,x)$ be
respectively the electric and magnetic fields. In dimensionless variables chosen so that the
speed of light, the charge and the mass of the particles are all equal to unity, the unknown
distribution function $f$ and electromagnetic field $(E,B)$ satisfy the Vlasov-Maxwell system
\begin{equation}
\label{RVM}
\begin{aligned}
\d_tf+v(\xi)\cdot\grad_xf&=-\Div_\xi[(E+v(\xi)\times B)f]\,,
\\
\d_tE-\curl_xB&=-j_f\,,\qquad\Div_xE=\rho_f\,,
\\
\d_tB+\curl_xE&=0\,,\,\quad\qquad\Div_xB=0\,.
\end{aligned}
\end{equation}
In this system, $\rho_f$ and $j_f$ denote respectively the charge and current densities
$$
\rho_f(t,x)=\int_{\R^3}f(t,x,\xi)d\xi\,,\quad
j_f(t,x)=\int_{\R^3}v(\xi)f(t,x,\xi)d\xi\,,
$$
while $v(\xi)$ is the relativistic velocity corresponding to a momentum $\xi$ for particles
with mass $1$ (in units such that the speed of light is $1$):
\begin{equation}
\label{v=}
v(\xi)=\frac{\xi}{\sqrt{1+|\xi|^2}}\,.
\end{equation}
The system (\ref{RVM}) is supplemented with the initial data
\begin{equation}
\label{IntlDt}
\begin{aligned}
f(0,x,\xi)&=f^{in}(x,\xi)\ge 0\,,&&\qquad x\in\R^3\,,\,\,\xi\in\R^3\,,
\\
(E,B)(0,x)&=(E^{in},B^{in})(x)\,,&&\qquad x\in\R^3\,.
\end{aligned}
\end{equation}

\subsection{Glassey-Strauss' conditional result}

In \cite{GS}, R. Glassey and W. Strauss have proved the following remarkable result: any
classical ($C^1$) solution of the system (\ref{RVM}) does not develop singularities as
long as the distribution function $f$ has compact support in the momentum variable $\xi$.
For each $f\equiv f(t,x,\xi)$, let
$$
R_f(t)=\inf\{r>0\,|\,f(t,x,\xi)=0\hbox{ for each }x\in\R^3\hbox{ and }|\xi|>r\}\,.
$$

\begin{Thm}[Glassey-Strauss \cite{GS}]\label{TH-GS}
Let $\tau>0$; let $f\in C^1([0,\tau)\times\R^3\times\R^3)$ and $(E,B)\in C^1([0,\tau)\times
\R^3)$ be a solution of (\ref{RVM}) with initial data $f^{in}\in C^1_c(\R^3\times\R^3)$ and
$(E^{in},B^{in})\in C^2_c(\R^3)$ satisfying the compatibility condition
\begin{equation}
\label{Cmpt}
\Div_xE^{in}=\int_{\R^3}f^{in}d\xi\,,\quad\Div_xB^{in}=0\,.
\end{equation}
If
$$
\varlimsup_{t\to\tau^-}\left(\|f(t)\|_{W^{1,\infty}_{x,\xi}}
    +\|(E,B)(t)\|_{W^{1,\infty}_x}\right)=+\infty
$$
then
$$
\varlimsup_{t\to\tau^-}R_f(t)=+\infty\,.
$$
\end{Thm}

The proof in \cite{GS} is based on a very ingenious argument: derivatives of the fields
$(E,B)$ with respect to $x$ are controlled in terms of $(\d_t+v(\xi)\cdot\grad_x)^mf$ for
$m=1,2$ and traded for derivatives of $f$ with respect to $\xi$ only. When computing the
charge and current densities, these $\xi$-derivatives disappear after integration by parts
in the variable $\xi$. However, this argument itself relies on rather formidable explicit
computations, especially for the derivatives of the fields $(E,B)$ in terms of $f$. In the
present paper, we give a shorter proof of Theorem \ref{TH-GS}. The simplifications come
mainly from
(a) expressing the field $(E,B)$ in terms of distributions of Lienard-Wiechert potentials
and
(b) a division lemma expressing second derivatives of the forward fundamental solution of
the wave equation in terms of the first and second power of the streaming operator acting
on that same fundamental solution.
In particular (b) is done in a direct and intrinsic way that avoids the repeated use of
Green's formula on the wave cone as in \cite{GS}. Moreover, (b) extends naturally to the
two-dimensional case --- which required a distinct treatment by the former method (see
\cite{GSc}).

As in \cite{BGPa}, the division lemma (Lemma \ref{LM-DVSN}) allows estimating the regularity
of $\xi$-averages of $u\equiv u(t,x,\xi)$ satisfying a coupled wave-transport system of the
form
$$
\Box_{t,x}u=f\,,\qquad(\d_t+v(\xi)\cdot\grad_x)f=P(t,x,\xi,D_\xi)g\,,
$$
where $g$ is given. The pseudo-differential approach in \cite{BGPa} leads to $L^p$ estimates
($1<p<+\infty$) in a quite general framework, being essentially based on the gap between the
characteristic speeds of $\Box_{t,x}$ and $(\d_t+v(\xi)\cdot\grad_x)$. The method described
in the present work uses in a deeper way the structure of the D'Alembert operator $\Box_{t,x}$,
as does \cite{GS}, and especially of its (forward) fundamental solution $Y$ in physical space
(i.e. in the $(t,x)$-variables). The decomposition (\ref{DrvtY2}) below plus the fact that $Y$
is a measure (in the case of space-dimension $3$) leads to $L^\infty$ estimates --- at the
only expense of a logarithmic term in $\grad_xf$ for the last piece $b_{ij}^2Y$ in that
decomposition: see subsection \ref{SBSCT-S3}. While these features were already present in
\cite{GS}, it seems that Lemma \ref{LM-DVSN} is new; in any case, its proof is based solely on
commutation properties of $\Box_{t,x}$ with the Lorentz boosts and not on the explicit form of
$Y$.

Another, new proof of Glassey-Strauss' conditional theorem has been recently given by
Klainerman-Staffilani \cite{KS}. Their proof is rather different from either that in
\cite{GS} or the present one. Maybe a combination of their arguments with the ones in
the present work could help in proving the global existence of classical solutions
without having to assume Glassey-Strauss' condition on $R_f$.

\section{Distributions of Lienard-Wiechert potentials}\label{LW-PT}

The formulation of (\ref{RVM}) involving distributions of Lienard-Wiechert potentials
appeared first in \cite{BGPa} and is recalled below. Let $u\equiv u(t,x,\xi)$ solve
$$
\Box_{t,x}u=f\,,\qquad u\rstr_{t=0}=\d_tu\rstr_{t=0}=0\,.
$$
Choose a vector field $A_I\equiv A_I(t,x)\in\R^3$ such that
\begin{equation}
\label{AI=}
\Div_xA_I=0\,,\quad \curl_xA_I=B^{in}\,,
\end{equation}
and solve for $A^0$ the wave equation
\begin{equation}
\label{A0=}
\Box A^0=0\,,\qquad A^0\rstr_{t=0}=A_I\,,\quad\d_tA^0\rstr_{t=0}=-E^{in}\,.
\end{equation}
Define the electromagnetic potential $(\phi,A)$ by
\begin{equation}
\label{Ptntl}
\phi(t,x)=\int_{\R^3}u(t,x,\xi)d\xi\,,\quad
    A(t,x)=A^0(t,x)+\int_{\R^3}v(\xi)u(t,x,\xi)d\xi\,,
\end{equation}
and the electromagnetic field $(E,B)$ by the usual formulas
\begin{equation}
\label{Flds}
E=-\d_tA-\grad_x\phi\,,\quad B=\curl_xA\,.
\end{equation}
Then the fields $(E,B)$ verify Maxwell's system of equations
$$
\begin{aligned}
\d_tE-\curl_xB&=-\int_{\R^3}v(\xi)fd\xi\,,&&\qquad\Div_xE=\int_{\R^3}fd\xi\,,
\\
\d_tB+\curl_xE&=0\,,&&\qquad\Div_xB=0\,,
\\
E\rstr_{t=0}&=E^{in}\,,&&\qquad B\rstr_{t=0}\,=B^{in}\,.
\end{aligned}
$$
Hence the relativistic Vlasov-Maxwell system (\ref{RVM}) can be put in the equivalent
form
\begin{equation}
\label{RVMLW}
\begin{aligned}
\d_tf+v(\xi)\cdot\grad_xf&=\Div_\xi[K_uf]\,,
\\
\Box_{t,x}u&=f\,,
\end{aligned}
\end{equation}
where $K_u$ is minus the Lorentz force field given by the formula
\begin{equation}
\label{K=}
\begin{aligned}
K_u(t,x,\xi)=\d_tA^0(t,x)&-v(\xi)\times\curl_xA^0(t,x)
\\
+\d_t\int_{\R^3}v(\xi)u(t,x,\xi)d\xi+\grad_x\int_{\R^3}ud\xi
    &-v(\xi)\times\curl_x\int_{\R^3}v(\xi)u(t,x,\xi)d\xi\,.
\end{aligned}
\end{equation}
with $A^0$ being defined by (\ref{A0=})-(\ref{AI=}). The initial conditions are
\begin{equation}
\label{IntlCdtnLW}
f\rstr_{t=0}=f^{in}\,,\quad u\rstr_{t=0}=0\,,\quad\d_tu\rstr_{t=0}=0\,,
\end{equation}
where $f^{in}$, $E^{in}$ and $B^{in}$ are assumed to satisfy the compatibility condition
(\ref{Cmpt}). Notice that the electromagnetic potential (\ref{Ptntl}) satisfies the Lorentz
gauge condition
\begin{equation}
\label{Lrnt}
\d_t\phi+\Div_xA=0\,.
\end{equation}
(Indeed, setting $G=\d_t\phi+\Div_xA$ and averaging Vlasov's equation in $\xi$, one has
$\Box G=\d_t\rho_f+\Div_xj_f=0$ with $G\rstr_{t=0}=\d_tG\rstr_{t=0}=0$).

\section{A division lemma}\label{DVSN-LM}

Let $Y\in\DD'(\R^4)$ be the forward fundamental solution of the d'Alember- tian, which is
characterized by
$$
\Box_{t,x}Y=\de_{(t,x)=(0,0)}\,,\qquad\hbox{supp}Y\subset\{(t,x)\in\R^4\,|\,|x|\le t\}\,.
$$
Although most of the present section can be understood without using the explicit formula
giving $Y$, we recall it below for convenience (with a slight abuse of notation):
\begin{equation}
\label{Y=}
Y(t,x)=\frac{\indc_{t>0}}{4\pi t}\,\delta(|x|-t)\,.
\end{equation}
Notice that the distribution $Y$ is homogeneous of degree $-2$ in $\R^4$. For $j=1,2,3$,
let $L_j=x_j\d_t+t\d_{x_j}$; one easily checks that $[\Box,L_j]=0$. On the other hand, one
has $L_j\de_{(t,x)=(0,0)}=0$, and this relation, together with the fact that $L_j$ commutes
with $\Box$, implies that
\begin{equation}
\label{LjY=}
L_jY=0\,,\qquad j=1,2,3\,.
\end{equation}
To each  $v\in\R^3$ is associated the streaming operator $T=\d_t+v\cdot\grad_x$. Let $\MM_m$
be the space of $C^\infty$ homogeneous functions of degree $m$ on $\R^4\setminus 0$. Below,
we use the notation
\begin{equation}
\label{x0nt}
x_0:=t\,,\qquad\hbox{ and }\d_j:=\d_{x_j}\,,\quad j=0,\ldots,3\,.
\end{equation}
The main result in the present section is

\begin{Lem}[Division lemma]\label{LM-DVSN}
For each $v\in\R^3$ such that $|v|<1$,

\begin{itemize}

\item there exists functions $a_i^k\equiv a_i^k(t,x)$ where $i=0,\ldots,3$ and $k=0,1$,
such that $a_i^k\in\MM_{-k}$ and
\begin{equation}
\label{DrvtY1}
\d_iY=T(a_i^0Y)+a_i^1Y\,,\qquad i=0,\ldots,3\,;
\end{equation}

\item there exists functions $b_{ij}^k\equiv b_{ij}^k(t,x)$ with $i,j=0,\ldots,3$, $k=0,1,2$,
such that $b_{ij}^k\in\MM_{-k}$ and
\begin{equation}
\label{DrvtY2}
\d_{ij}^2Y=T^2(b_{ij}^0Y)+T(b_{ij}^1Y)+b_{ij}^2Y\,,\qquad i,j=0,\ldots,3\,;
\end{equation}
\item moreover, the functions $b_{ij}^2$ satisfy the conditions
\begin{equation}
\label{Av=0}
\int_{\S^2}b_{ij}^2(1,y)d\si(y)=0\,,\qquad i,j=0,\ldots,3\,,
\end{equation}
where $d\si(y)$ is the rotation invariant surface element on the unit sphere $\S^2$ of $\R^3$. In
both formulas (\ref{DrvtY1}) and (\ref{DrvtY2}), $a^0_iY$, $a^1_iY$, $b^0_{ij}Y$ and $b^1_{ij}Y$
designate, for each $i,j=0,\ldots,3$, the unique extensions\footnote{We abandon in the main body
of the text the notation $\dot{f}$ for the unique homogeneous extension to $\R^4$ of a distribution
$f$ on $\R^4\setminus 0$ that is homogeneous of degree $>-4$; this notation is used in the appendix
only for the sake of clarity.} as homogeneous distributions on $\R^4$ of those same expressions ---
which are {\it a priori} only defined on $\R^4\setminus 0$. Likewise, $b_{ij}^2Y$ designates, for
$i,j=0,\ldots,3$ the unique extension as a homogeneous distribution of degree $-4$ on $\R^4$ of that
same expressions for which the relation (\ref{DrvtY2}) holds in the sense of distributions on $\R^4$.
\end{itemize}
\end{Lem}

\smallskip
\noindent
{\bf Remark.} The first and second statements in Lemma \ref{LM-DVSN} hold verbatim in space dimension
2. As for the third statement, the degree of homogeneity of $b^2_{ij}Y$ is $-3$ in $\R^3$ in the case
of space dimension 2, and the condition (\ref{Av=0}) becomes
\begin{equation}
\label{Av02D}
\int_{|y|<1}\frac{b_{ij}^2(1,y)}{\sqrt{1-|y|^2}}\,dy=0\,,\qquad i,j=0,\ldots,2\,.
\end{equation}

\begin{proof}
Observe that
$$
\begin{aligned}
\sum_{j=1}^3v_jL_j&=\,\,(x\cdot v-t)\d_t+tT\,,
\\
(t-x\cdot v)L_i+x_i\sum_{j=1}^3v_jL_j&=t[(t-x\cdot v)\d_i+x_iT]\,,\quad i=1,2,3\,.
\end{aligned}
$$
These relations and (\ref{LjY=}) imply that
\begin{equation}
\label{dYTY}
(t-x\cdot v)\d_tY=tTY\,,\quad t(x\cdot v-t)\d_iY=tx_iTY\,,\quad i=1,2,3\,.
\end{equation}
Set
$$
\a_0(t,x)=\frac{t}{t-x\cdot v}\,,\qquad \a_i(t,x)=\frac{x_i}{x\cdot v-t}\,,\,\,i=1,2,3.
$$
Since $|v|<1$, the functions $\a_i$ are $C^\infty$ near the support of the restriction
of $Y$ to $\R^4\setminus 0$: hence $\a_iTY$ defines, for $i=0,\ldots,3$, a distribution
on $\R^4\setminus 0$ that is homogeneous of degree $-3$. It has a unique extension as a
homogeneous distribution of degree $-3$ on $\R^4$, still denoted by $\a_iTY$. Because
of (\ref{dYTY}), the distribution $\d_iY-\a_iTY$ has support in the set $\{(t,x)\in\R^4
\,|\,x\cdot v=t\}\cup\{(t,x)\in\R^4\,|\,t=0\}$; since $Y$ is supported in the wave cone
$\{(t,x)\in\R^4\,|\,|x|\le t\}$, it follows that
$$
\begin{aligned}
\Supp(\d_iY-\a_iTY)\subset&\{(t,x)\in\R^4 \,|\,x\cdot v=t\hbox{ or }t=0\}
\\
&\cap\{(t,x)\in\R^4\,|\,|x|\le t\}=\{(0,0)\}\,.
\end{aligned}
$$
Thus $\d_iY-\a_iTY$ is both a homogeneous distribution on $\R^4$ of degree $-3$ and
a finite linear combination of $\delta_{(t,x)=(0,0)}$ and of its derivatives: hence
$$
\d_iY-\a_iTY=0\,,\quad i=0,\ldots,3\,.
$$
The same holds if one replaces $\a_i$ by a smooth truncation $a_i^0$ of it near its
singular set: indeed, as observed above, this singular set $\{(t,x)\in\R^4 \,|\,x
\cdot v=t\}$ does not intersect the support of $Y$ restricted to $\R^4\setminus 0$.
Hence $\d_iY=a_i^0TY=T(a_i^0Y)-T(a_i^0)Y$ and formula (\ref{DrvtY1}) holds with
\begin{equation}
\label{dfnt-a}
a_i^0(t,x)=\a_i(t,x)\chi\left(\frac{|x|}t\right)\,,\quad\hbox{ and }a_i^1=-Ta_i^0\,,
    \quad i=0,\ldots,3\,,
\end{equation}
where $\chi\in C^\infty_c(\R_+)$ satisfies
$$
0\le\chi\le 1\,,\quad \chi\rstr_{\left[0,\tfrac12+\tfrac1{2|v|}\right]}\equiv 1\,,
    \quad\hbox{supp}\chi\subset\left[0,\tfrac1{|v|}\right[\,.
$$
By the same argument, the equality
\begin{equation}
\label{DrvtmY}
\d_i(mY)=T(ma_i^0Y)+\left(\d_im-T(ma_i^0)\right)Y\,,\qquad i=0,\ldots,3
\end{equation}
holds in the sense of distributions on $\R^4$ for each $m\in\MM_0$ --- where $mY$,
$ma_i^0Y$ and $(\d_im-T(ma_i^0))Y$ designate the homogeneous extensions to $\R^4$ of these
same distributions that are defined and homogeneous of degree $>-4$ on $\R^4\setminus 0$.

If $m\in\MM_{-1}$, the distributions $mY$ and $ma_i^0Y$ for $i=0,\ldots,3$ are homogeneous
of degree $-3$ in $\R^4\setminus 0$ and thus have unique extensions to $\R^4$ as homogeneous
distributions of degree $-3$ (see appendix). Since
$$
\d_i(mY)-T(ma_i^0Y)=\left(\d_im-T(ma_i^0)\right)Y\,,\qquad i=0,\ldots,3
$$
in the sense of distributions on $\R^4\setminus 0$, the r.h.s. of the above equality extends
as a homogeneous distribution of degree $-4$ on $\R^4$. Hence (see appendix)
$$
\hbox{Res}_0\left(\d_im-T(ma_i^0)\right)Y=0\,.
$$
Using (\ref{RsdInt}) and the formula for $Y$ in the case of space dimension 3, we can write
this condition in the following, more explicit form: for each $\chi\in C_c^\infty((0,+\infty))$
$$
\begin{aligned}
{}&\int_{\R^4}\left(\d_im-T(ma_i^0)\right)(t,x)\chi(t^2+|x|^2)Y(t,dx)dt
\\
=
&\int_0^{+\infty}\frac{\chi(2t^2)}{4\pi t}
    \int_{|x|=t}\left(\d_im-T(ma_i^0)\right)(t,x)d\si_t(x)dt\,.
=0
\end{aligned}
$$
Here, $d\si_t(x)$ designates the surface element on the sphere of equation $|x|=t$. Also, in the
case of space dimension 3, the distribution $Y$ is in fact a measure --- we recall from (\ref{Y=})
that $Y(t,dx)\equiv\frac{\indc_{t>0}}{4\pi t}\,d\si_t(x)$ --- which makes it legitimate to write
the left-hand side of the equality above as an integral.

The function $\d_im-T(ma_i^0)$ is homogeneous of degree $-2$ on $\R^4\setminus 0$, so that,
in terms of the new variable $y=x/t$, the last integral in the relation above becomes
$$
\begin{aligned}
{}&\int_0^{+\infty}\frac{\chi(2t^2)}{4\pi t}
    \int_{|x|=t}\left(\d_im-T(ma_i^0)\right)(t,x)d\si_t(x)dt
\\
=
&\int_0^{+\infty}\frac{\chi(2t^2)}{4\pi t}dt
    \int_{|y|=1}\left(\d_im-T(ma_i^0)\right)(1,y)d\si(y)\,.
\end{aligned}
$$
Eventually
\begin{equation}
\label{Avm0}
\int_{\S^2}\left(\d_im-T(ma_i^0)\right)(1,x)d\si(x)=0\,.
\end{equation}

In order to obtain (\ref{DrvtY2}), we apply $\d_j$ to both sides of the relation
(\ref{DrvtY1}) so as to obtain
$$
\d_{ij}^2Y=T(\d_i(a_j^0Y)) +\d_i(a_j^1Y)\,,\quad i,j=0,\ldots,3\,.
$$
Then we apply (\ref{DrvtmY}), first with $m=a_j^0\in\MM_0$, then with $m=a_j^1=-Ta_j^0
\in\MM_{-1}$: this leads (\ref{DrvtY1}) to (\ref{DrvtY2}) with
\begin{equation}
\label{dfnt-b}
\begin{aligned}
b_{ij}^0&=a_i^0a_j^0\,,
\\
b_{ij}^1&=\d_ia_j^0-T(a_i^0a_j^0)+a_i^0a_j^1\,,\quad i,j=0,\dots,3\,,
\\
b_{ij}^2&=\d_ia_j^1-T(a_i^0a_j^1)\,,
\end{aligned}
\end{equation}
the functions $a_i^k$ being defined in (\ref{dfnt-a}).

Finally, the condition (\ref{Avm0}) with $m=a_j^1$ is equivalent to (\ref{Av=0}).
\end{proof}

\section{Bounds on the electro-magnetic field}\label{BD-FLD}

After these preparations, we give the proof of Theorem \ref{TH-GS}. Since the solution
$(f,E,B)$ belongs to $C^1([0,\tau)\times\R^3\times\R^3)\times C^1([0,\tau)\times\R^3)^2$,
the distribution function $f$ is constant along the characteristic curves of the vector
field $(v(\xi),-K_u(t,x,\xi))$ --- observe that $\Div_\xi K_u(t,x,\xi)=0$. In particular
\begin{equation}
\label{MxPr-f}
\|f\|_{L^\infty([0,\tau]\times\R^3\times\R^3)}=\|f^{in}\|_{L^\infty(\R^3\times\R^3)}\,,
\end{equation}
and, since $f^{in}$ has compact support,
$$
\sup_{t\in[0,\tau']}R_f(t)<+\infty\hbox{ for each }\tau'\in[0,\tau)\,.
$$
Hence proving Theorem \ref{TH-GS} amounts to proving the implication
\begin{equation}
\label{ImplGS}
\sup_{t\in [0,\tau)}R_f(t)<+\infty\,\Longrightarrow
    \sup_{t\in[0,\tau)}\left(\|f(t)\|_{W^{1,\infty}_{x,\xi}}
    +\|(E,B)(t)\|_{W^{1,\infty}_x}\right)<+\infty\,.
\end{equation}

From now on, assume that
\begin{equation}
\label{AssmR}
\sup_{t\in [0,\tau)}R_f(t)=r^*\,.
\end{equation}
In other words
\begin{equation}
\label{Supp-xi}
f(t,x,\xi)\equiv 0\hbox{ and }u(t,x,\xi)\equiv 0\,,
    \quad t\in[0,\tau)\,,\,\,x\in\R^3\,,\,\,|\xi|>r^*\,.
\end{equation}

Next we want to estimate the electromagnetic field in $L^\infty([0,\tau)\times\R^3)$.
Start from the relation\footnote{In the sequel, the notation $f\star g$ always means
convolution in $\R^4_{t,x}$; the symbol $\star_x$ designates the convolution in the
variable $x\in\R^3$ only.} $u=Y\star(\indc_{t\ge 0}f)$ and use Lemma \ref{LM-DVSN} to
compute, for each $m\equiv m(\xi)$ in $C(\R^3)$
$$
\begin{aligned}
\d_j\int m(\xi)u(t,x,\xi)d\xi
    &=\int m(\xi)\left(\d_jY\star(\indc_{t\ge 0}f)\right)(t,x,\xi)d\xi
\\
&=\int m(\xi)\left((a_j^0Y)\star T(\indc_{t\ge 0}f)\right)(t,x,\xi)d\xi
\\
&+\int m(\xi)\left((a_j^1Y)\star(\indc_{t\ge 0}f)\right)(t,x,\xi)d\xi
\end{aligned}
$$
for $j=0,\ldots,3$, with $a_j^k\equiv a_j^k(t,x,\xi)$ given by (\ref{dfnt-a}) for  $v
\equiv v(\xi)$ as in (\ref{v=}). First, $a_j^k\in C^\infty((\R^4\setminus 0)\times\R^3)$;
also $\d_\xi^\beta a_j^k(\cdot,\cdot,\xi)$ is an element of $\MM_{-k}$ for each $\xi\in
\R^3$ and each multi-index $\beta\in\N^3$. By the first equation in (\ref{RVMLW}),
$$
T(\indc_{t\ge 0}f)=\de_{t=0}f^{in}+\indc_{t\ge 0}\Div_\xi(K_uf)\,;
$$
hence, if $m\in W^{1,\infty}$, one finds that
$$
\begin{aligned}
\int m(\xi)&\left((a_j^0Y)\star T(\indc_{t\ge 0}f)\right)(t,x,\xi)d\xi
\\
&=
\int\left(\left(-\grad_\xi(ma_j^0)Y\right)\star(\indc_{t\ge 0}K_uf)\right)(t,x,\xi)d\xi
\\
&+\int m(\xi)
    \left(\left(a_j^0(t,\cdot,\cdot)Y(t,\cdot)\right)\star_x f^{in}\right)(x,\xi)d\xi\,.
\end{aligned}
$$
Let $\phi\in C_c^\infty(\R^3)$ satisfy
\begin{equation}
\label{CndTrnc}
\phi\ge 0\,,\quad \phi(\xi)=1\hbox{ for }|\xi|\le r^*\,,\quad\phi(\xi)=0
    \hbox{ for }|\xi|\ge 2r^*\,.
\end{equation}
By the support condition (\ref{Supp-xi}), for each $m\in C(\R^3)$, one has
\begin{equation}
\label{Lclz}
\begin{aligned}
\int_{\R^3}m(\xi)f(t,x,\xi)d\xi&=\int_{\R^3}\phi(\xi)m(\xi)f(t,x,\xi)d\xi\,,
\\
\int_{\R^3}m(\xi)u(t,x,\xi)d\xi&=\int_{\R^3}\phi(\xi)m(\xi)u(t,x,\xi)d\xi\,.
\end{aligned}
\end{equation}
Since $Y(t,\cdot)$ is a positive measure with total mass $t$, it follows from (\ref{Lclz})
that
\begin{equation}
\label{djSmu}
\begin{aligned}
{}&\left|\d_j\int m(\xi)u(t,x,\xi)d\xi\right|
\\
&\le
\|m\|_{W^{1,\infty}}\|\phi a_j^0\|_{L^\infty_{t,x}(W^{1,\infty}_\xi)}\tfrac43 \pi r^{*3}
    \int_0^t(t-s)\|fK_u(s,\cdot,\cdot)\|_{L^\infty}ds
\\
&+
\|m\|_{L^\infty}\|\phi ta_j^1\|_{L^\infty}\tfrac43 \pi r^{*3}
    \int_0^t\|f(s,\cdot,\cdot)\|_{L^\infty}ds
\\
&+
\|m\|_{L^\infty}\|\phi a_j^0\|_{L^\infty}\tfrac43 \pi r^{*3}t\|f^{in}\|_{L^\infty}\,.
\end{aligned}
\end{equation}
Without loss of generality, we only consider the case where $B^{in}=0$; hence $A^0=-
Y(t,\cdot)\star_xE^{in}\in C_t(W^{2,\infty}_x)$. We recall at this point the elementary
estimates that hold for $k=0,1,2$:
\begin{equation}
\label{EstmA0}
\begin{aligned}
\|A^0(t)\|_{W^{k,\infty}_x}&\le t\|E^{in}\|_{W^{k,\infty}}
\\
\|\d_tA^{0}(t)\|_{W^{k-1,\infty}_x}&\le (1+t)\|E^{in}\|_{W^{k,\infty}}\,.
\end{aligned}
\end{equation}
Denote
$$
I_m(t)=\sup_{j=0,\ldots,3}\left\|\d_j\int m(\xi)u(t,\cdot,\xi)d\xi\right\|_{L^\infty}\,;
$$
by using (\ref{EstmA0}), (\ref{djSmu}) and (\ref{K=}), one sees, for each $t\in[0,\tau)$
and some positive constant $C(\tau,r^*,\|m\|_{W^{1,\infty}},\|f^{in}\|_{L^\infty})>0$,
that
\begin{equation}
\label{Grnw1}
I_m(t)\le C(\tau,r^*,\|m\|_{W^{1,\infty}},\|f^{in}\|_{L^\infty})
    \left(1+\int_0^t\left(I_1(s)+I_v(s)\right)ds\right)\,.
\end{equation}
Using (\ref{Grnw1}) for $m\equiv 1$ and $m=v$ and applying Gronwall's inequality, one
finds that
\begin{equation}
\label{BndFld}
\sup_{t\in[0,\tau)}I_m(t)<+\infty\quad\hbox{ for each }m\in W^{1,\infty}(\R^3)\,.
\end{equation}
In particular, using again (\ref{K=}), one finds eventually that
\begin{equation}
\label{BndK}
\|K_u\|_{L^\infty([0,\tau)\times\R^3;W^{k,\infty}_\xi)}<+\infty\hbox{ for each }k\ge 0\,.
\end{equation}

\section{Bounds on first derivatives}\label{BD-DRVT}

For each $m\in C(\R^3)$, one has, by using Lemma \ref{LM-DVSN} and the support condition
(\ref{Supp-xi}),
$$
\begin{aligned}
\d_{ij}
\int m(\xi)u(t,x,\xi)d\xi&=\int m(\xi)\d_{ij}Y\star(\indc_{t\ge 0}f)(t,x,\xi)d\xi
\\
&=\int m(\xi)\left((b_{ij}^0Y)\star T^2(\indc_{t\ge 0}f)\right)(t,x,\xi)d\xi
\\
&+\int m(\xi)\left((b_{ij}^1Y)\star T(\indc_{t\ge 0}f)\right)(t,x,\xi)d\xi
\\
&+\int m(\xi)\left((b_{ij}^2Y)\star(\indc_{t\ge 0}f)\right)(t,x,\xi)d\xi
\\
&=S_1+S_2+S_3\,.
\end{aligned}
$$
for $j=0,\ldots,3$, where $b_{ij}^k\equiv b_{ij}^k(t,x,\xi)$ is given by (\ref{dfnt-b}). In
the second integral appearing in the r.h.s. of the relation above, $T(\indc_{t\ge 0}f)$ is
replaced by $\de_{t=0}f^{in}+\indc_{t\ge 0}\Div_\xi(K_uf)$, in view of the first equation
in (\ref{RVMLW}) in the previous section. Likewise, in the first integral in the r.h.s. of
the equality above, $T^2(\indc_{t\ge 0}f)$ is expressed as
$$
\begin{aligned}
T^2(\indc_{t\ge 0}f)&=T(\de_{t=0}f^{in})+T\left(\indc_{t\ge 0}\Div_\xi(K_uf)\right)
\\
&=\de'_{t=0}f^{in}+\de_{t=0}\left(v\cdot\grad_xf^{in}+\Div_\xi(K_u^{in}f^{in})\right)
\\
&+\indc_{t\ge 0}\Div_\xi\left(fTK_u+K_u\Div_\xi(K_uf)\right)
    +\indc_{t\ge 0}[T,\Div_\xi](K_uf)
\\
&=\de'_{t=0}f^{in}+\de_{t=0}\left(v\cdot\grad_xf^{in}+\Div_\xi(K_u^{in}f^{in})\right)
\\
&+\indc_{t\ge 0}\grad_\xi^{\otimes 2}:(fK_u^{\otimes 2})
    +\indc_{t\ge 0}\Div_\xi(fTK_u-fK_u\cdot\grad_\xi K_u)
\\
&-\left(\grad_\xi v\right)^T:\grad_x\left(\indc_{t\ge 0}fK_u\right)\,.
\end{aligned}
$$
Below, we shall use the notation
$$
J_m(t)=\sup_{i,j=0,\ldots,3}\left\|\d_{ij}\int m(\xi)u(t,\cdot,\xi)d\xi\right\|_{L^\infty}\,.
$$

\subsection{Estimating $S_1$}

We decompose further
$$
\begin{aligned}
S_1&=\int m(\xi)(b_{ij}^0Y)\star \left(\de'_{t=0}f^{in}
    +\de_{t=0}(v\cdot\grad_xf^{in}+\Div_\xi(K_u^{in}f^{in}))\right)d\xi
\\
+&
\int\left(\left(\grad_\xi^{\otimes 2}(mb_{ij}^0)Y\right)\star
    \left(\indc_{t\ge 0}fK_u^{\otimes 2}\right)\right)(t,x,\xi)d\xi
\\
+&
\int\left(\left(-\grad_\xi(mb_{ij}^0)Y\right)\star
    \left(\indc_{t\ge 0}(fTK_u-fK_u\cdot\grad_\xi K_u)\right)\right)(t,x,\xi)d\xi
\\
+&
\int m(\xi)\left(\left(\grad_\xi v\cdot\grad_x(b_{ij}^0Y)\right)\star
    \left(\indc_{t\ge 0}fK_u\right)\right)(t,x,\xi)d\xi
\\
=&
S_{11}+S_{12}+S_{13}+S_{14}\,.
\end{aligned}
$$
By using the classical estimates (\ref{EstmA0}) for the wave equation, together with
the support condition (\ref{Supp-xi}) and the definition of $\phi$ in (\ref{CndTrnc})
\begin{equation}
\label{S11}
\begin{aligned}
|S_{11}|&\le\|\phi mb_{ij}^0\|_{L^\infty_x(W^{1,\infty}_{t,\xi})}
\\
&\times\tfrac43\pi r^{*3}(1+\tau)^2\|f^{in}\|_{W^{1,\infty}}
\left(1+\|K_u^{in}\|_{L^{\infty}([0,\tau)\times\R^3;W^{1,\infty}_\xi)}\right)\,.
\end{aligned}
\end{equation}
By the same argument as in section \ref{BD-FLD}
\begin{equation}
\label{S12}
|S_{12}|\le
\|\phi mb_{ij}^0\|_{L^\infty_{t,x}(W^{2,\infty}_\xi)}\tfrac43\pi r^{*3}
\tfrac12\tau^2\|f^{in}\|_{L^\infty}\|K_u\|^2_{L^\infty([0,\tau)\times\R^3\times\R^3)}\,.
\end{equation}
Likewise
\begin{equation}
\label{S13}
\begin{aligned}
|S_{13}|&\le\|\phi mb_{ij}^0\|_{L^\infty_{t,x}(W^{1,\infty}_\xi)}
    \tfrac43\pi r^{*3}\|f^{in}\|_{L^\infty}
\\
&\times\left(\int_0^t(t-s)(J_1(s)+J_v(s))ds
    +\|K\|^2_{L^\infty([0,\tau)\times\R^3;W^{1,\infty}_\xi)}\right)\,.
\end{aligned}
\end{equation}
In $S_{14}$, we apply once more Lemma \ref{LM-DVSN} --- or (\ref{DrvtmY}) with $m=b_{ij}^0
\in\MM_0$ and $j=1,2,3$ --- so as to write
$$
\d_k(b_{ij}^0Y)=T(b_{ij}^0a_k^0Y)-\left(\d_kb_{ij}^0-T(b_{ij}^0a_k^0)\right)Y\,;
$$
replacing this in the expression giving $S_{14}$ and proceeding as in section \ref{BD-FLD}
leads to
\begin{equation}
\label{S14}
\begin{aligned}
|S_{14}|&\le\tfrac43\pi r^{*3}\|f^{in}\|_{L^\infty}\left[\sup_{k=0,\ldots,3}
    \|\phi ma_k^0b_{ij}^0\|_{L^\infty_{t,x}(W^{1,\infty}_\xi)}\right.
\\
&\times
\left(1+\|K_u\|_{L^\infty([0,\tau)\times\R^3;W^{1,\infty}_\xi)}\right)
    \int_0^t(t-s)(J_1(s)+J_v(s))ds
\\
&+
\sup_{k=0,\ldots,3}
\left\|\phi m\left(t\d_kb_{ij}^0-tT(b_{ij}^0a_k^0)\right)\right\|_{L^\infty}
    \|K_u\|_{L^\infty([0,\tau)\times\R^3\times\R^3)}
\\
&+\left.
\sup_{k=0,\ldots,3}
\|\phi ma_k^0b_{ij}^0\|_{L^\infty}\tfrac12t^2\|K^{in}\|_{L^\infty}\right]\,.
\end{aligned}
\end{equation}

\subsection{Estimating $S_2$}

This part of the argument follows section \ref{BD-FLD}, except $b_{ij}^1\in\MM_{-1}$ while
$a_{i}^0\in\MM_0$. Thus
$$
\begin{aligned}
S_2&=
\int\left(\left(-\grad_\xi(mb_{ij}^1)Y\right)\star (\indc_{t\ge 0}K_uf)\right)(t,x,\xi)d\xi
\\
&+\int m(\xi)\left((b_{ij}^1(t,\cdot,\cdot)Y(t,\cdot))\star_xf^{in}\right)(x,\xi)d\xi
\end{aligned}
$$
so that, by the same estimates leading to the last two terms in the r.h.s. of (\ref{djSmu}),
one arrives at
\begin{equation}
\label{S2}
|S_2|\le
\|\phi mtb_{ij}^1\|_{L^\infty_{t,x}(W^{1,\infty}_\xi)}\tfrac43\pi r^{*3}
\left(\tau\|K_u\|_{L^\infty([0,\tau)\times\R^3\times\R^3)}+\|f^{in}\|_{L^\infty}\right)\,.
\end{equation}

\subsection{Estimating $S_3$}\label{SBSCT-S3}

Let $\phi\in C^\infty_c(\R^4\setminus 0)$; since $b^2_{ij}$ is homogeneous of degree $-2$
(we recall from Lemma \ref{LM-DVSN} that $b_{ij}^2(\cdot,\cdot,\xi)\in\MM_{-2}$), one has
$$
\begin{aligned}
\la b^2_{ij}Y,\phi\ra&=\int_0^{+\infty}
    \int_{\S^2}\frac1{4\pi t}b^2_{ij}(t,t\om,\xi)\phi(t,t\om)t^2d\si(\om) dt
\\
&=
\int_0^{+\infty}\int_{\S^2}\frac1{4\pi t}b^2_{ij}(1,\om,\xi)\phi(t,t\om)d\si(\om)dt\,.
\end{aligned}
$$
(where $d\si(\om)$ is the rotation invariant surface element on $\S^2$). Further, the relation
(\ref{Av=0}) shows that, for each $\psi\in C^\infty_c(\R^4)$, the quantity
\begin{equation}
\label{VP}
\begin{aligned}
\la\VP(b^2_{ij}Y),\psi\ra
&=\int_\th^{+\infty}\int_{\S^2}b^2_{ij}(1,\om,\xi)\frac{\psi(t,t\om)}{4\pi t}d\si(\om)dt
\\
&+
\int_0^{\th}\int_{\S^2}b^2_{ij}(1,\om,\xi)\frac{\psi(t,t\om)-\psi(t,0)}{4\pi t}d\si(\om)dt
\end{aligned}
\end{equation}
is independent of $\th\in\R_+$. This defines $\VP(b^2_{ij}Y)$ as a homogeneous distribution
of degree $-4$ on $\R^4$ that extends $b^2_{ij}Y\rstr_{\R^4\setminus 0}$. Hence (see appendix)
\begin{equation}
\label{b-VPb}
b^2_{ij}(\cdot,\cdot,\xi)Y-\VP\left(b^2_{ij}(\cdot,\cdot,\xi)Y\right)
    =c_{ij}(\xi)\de_{(t,x)=(0,0)}\,,
\end{equation}
where $c_{ij}\in C^\infty(\R^3)$ --- we recall that the l.h.s. of the equality above is of
class $C^\infty$ in $\xi$.

Therefore, for $\th_t\in(0,t)$ to be chosen later,
$$
\begin{aligned}
S_3&\!-\!\int m(\xi)c_{ij}(\xi)f(t,x,\xi)d\xi
\\
&=\!\int m(\xi)\left(\VP(b_{ij}^2Y)\star(\indc_{t\ge 0}f)\right)(t,x,\xi)d\xi
\\
&=\!\int\!\!m(\xi)\!\!
\int_{\th_t}^t\!\int_{\S^2}b^2_{ij}(1,\om,\xi)f(t-s,x-s\om,\xi)\frac{d\si(\om)ds}{4\pi s}d\xi
\\
&+\!\int\!\!m(\xi)\!\!
\int_0^{\th_t}\!\!\!\int_{\S^2}b^2_{ij}(1,\om,\xi)
    \frac{f(t\!-\!s,x\!-\!s\om,\xi)\!-\!f(t\!-\!s,x,\xi)}{4\pi s}d\si(\om)ds
    d\xi.
\end{aligned}
$$
The first integral in the r.h.s. of the relation above is estimated by
$$
\begin{aligned}
\left|\int_{\th_t}^t\!\int_{\S^2}
    b^2_{ij}(1,\om,\xi)f(t-s,x-s\om,\xi)\frac{d\si(\om)ds}{4\pi s}\right|&
\\
\le
\ln(t/\th_t)\|b_{ij}^2(1,\cdot,\xi)\|_{L^\infty(\S^2)}\|f\|_{L^\infty}&\,,
\end{aligned}
$$
while the second is estimated by
$$
\begin{aligned}
\left|\int_0^{\th_t}\!\!\int_{\S^2}b^2_{ij}(1,\om,\xi)
    \frac{f(t-s,x-s\om,\xi)-f(t-s,x,\xi)}{4\pi s}d\si(\om)ds\right|&{}
\\
\le\th_t\|b_{ij}^2(1,\cdot,\xi)\|_{L^\infty(\S^2)}
    \|\grad_xf\|_{L^\infty([0,t]\times\R^3\times\R^3)}&\,.
\end{aligned}
$$
Choosing
$$
\th_t=\inf\left(\frac1{\|\grad_xf\|_{L^\infty([0,t]\times\R^3\times\R^3)}},t\right)
$$
one finds that
\begin{equation}
\label{S3}
\begin{aligned}
|S_3|&\le Cr^{*3}\|m\|_{L^\infty}\left[\|c_{ij}\|_{L^\infty(B(0,r^*))}\|f\|_{L^\infty}
+\|b_{ij}^2(1,\cdot,\cdot)\|_{L^\infty(\S^2\times\R^3)}\right.
\\
&\left.\times\left(1+\|f\|_{L^\infty}
\ln_+\left(t\|\grad_xf\|_{L^\infty([0,t]\times\R^3\times\R^3)}\right)\right)\right]
\end{aligned}
\end{equation}
where $\ln_+z=\sup(\ln z,0)$.

\smallskip
\noindent
{\bf Remark.} In the case of space dimension 2, a similar argument, based on (\ref{Av02D})
instead of (\ref{Av=0}), leads to an estimate that involves a logarithmic term just as in
(\ref{S3}). The condition (\ref{Av02D}) is not apparent in \cite{GSc}, which uses instead
the fact that $b^2_{ij}Y$ is a linear combination of derivatives of distributions that are
homogeneous of degree $\ge -2$ in $\R^3$ --- see p. 344 of \cite{GSc}. As explained in the
appendix, this is equivalent to (\ref{Av02D}).

\subsection{Proof of Theorem \ref{TH-GS}}

The estimates (\ref{S11}), (\ref{S12}), (\ref{S13}), (\ref{S14}), (\ref{S2}) and (\ref{S3})
show that, for each $m\in W^{2,\infty}(\R^3)$, there exists a positive constant $C_2\equiv
C_2(\tau,r^*,\|m\|_{W^{2,\infty}},\|f^{in}\|_{W^{1,\infty}})$ such that
\begin{equation}
\label{Jm}
\begin{aligned}
J_m(t')&\le C_2(\tau,r^*,\|m\|_{W^{2,\infty}},\|f^{in}\|_{W^{1,\infty}})
\\
&\times\left(1+\int_0^{t'}(J_1(s)+J_m(s))ds
+\ln_+\left(\|\grad_xf\|_{L^\infty([0,t]\times\R^3\times\R^3)}\right)\right)
\end{aligned}
\end{equation}
for each $t$ and $t'$ such that $0<t'<t<\tau$. Using (\ref{Jm}) with $m\equiv 1$ and $m=v$
and applying Gronwall's inequality shows that, for each $t\in[0,\tau)$
\begin{equation}
\label{J1v}
J_1(t)+J_v(t)\le 2C_2e^{2C_2\tau}\left(1+
\ln_+\left(\|\grad_xf\|_{L^\infty([0,t]\times\R^3\times\R^3)}\right)\right)\,.
\end{equation}
In particular, this implies the existence of another positive constant $C_3\equiv C_3(\tau,
r^*,\|m\|_{W^{2,\infty}},\|f^{in}\|_{W^{1,\infty}})$ such that
\begin{equation}
\label{Grnw2}
\|K_u(t)\|_{W^{1,\infty}_{x,\xi}}\le C_3e^{2C_2\tau}\left(1+
\ln_+\left(\|\grad_xf\|_{L^\infty([0,t]\times\R^3\times\R^3)}\right)\right)\,.
\end{equation}
Finally, differentiating in $(x,v)$ the transport equation in (\ref{RVMLW}) and integrating
in $t$ shows that
\begin{equation}
\label{EstmDf}
\begin{aligned}
\|\grad_{x,\xi}f(t)\|_{L^\infty_{x,\xi}}
&\le
\|\grad_{x,\xi}f^{in}\|_{L^\infty_{x,\xi}}
\\
&+
\int_0^t\left(\|\grad_\xi v\|_{L^\infty}+\|\grad_{x,\xi}K_u(s)\|_{L^{\infty}_{x,\xi}}\right)
\|\grad_{x,\xi}f(s)\|_{L^\infty_{x,\xi}}ds\,.
\end{aligned}
\end{equation}
The estimates (\ref{Grnw2}) and (\ref{EstmDf}) show that the Lipschitz semi-norm of $f$, i.e.
$N(t)=\sup_{s\in[0,t]}\|\grad_{x,\xi}f(s)\|_{L^\infty_{x,\xi}}$ satisfies a logarithmic
Gronwall inequality of the form
$$
N(t)\le N(0)+C\int_0^t(1+\ln_+N(s))N(s)ds\,,\quad t\in[0,\tau]\,.
$$
This implies that $N\in L^\infty([0,\tau])$. Inserting this in (\ref{J1v}) and using
(\ref{Flds}) shows that $(E,B)\in L^\infty([0,\tau],W^{1,\infty}(\R^3))$, which in
turn implies Theorem \ref{TH-GS}.

\section{Appendix: Homogeneous distributions}

This section recalls some classical material from \cite{Ge} (chapter III, section 3.3) and
\cite{Ho} (pp. 75--79 and Theorem 3.2.3).

A distribution $f$ on $\R^N$ (resp. $\R^N\setminus 0$) is homogeneous of degree $\a$ if
$\la f,M_\l\phi\ra=\l^{\a+N}\la f,\phi\ra$ (where $M_\l\phi(x)=\phi(x/\l)$) for each
$\l>0$ and each $\phi\in C^\infty_c(\R^N)$ (resp. $\phi\in C^\infty_c(\R^N\setminus 0)$).
Equivalently, $f\in\DD'(\R^N)$ (resp. of $\DD'(\R^N\setminus 0)$) is homogeneous of
degree $\a$ if and only if
\begin{equation}
\label{Euler}
\Div_x(xf)=(\a+N)f\hbox{ on }\R^N\hbox{ (resp. on }\R^N\setminus 0)
\end{equation}
in the sense of distributions (Euler's relation in conservation form).

For $\a>-N$, each homogeneous distribution $f$ of degree $\a$ on $\R^N\setminus 0$ has a
unique extension $\dot{f}$ that is a homogeneous distribution on $\R^N$.

If $f\in\DD'(\R^N\setminus 0)$ is homogeneous of degree $-N$, $\Div_x((xf)^\cdot)$ is a
homogeneous distribution of degree $-N$ on $\R^N$ supported in $\{0\}$, hence there exists
$c\in\R$
\begin{equation}
\label{Rsd}
\Div_x((xf)^\cdot)=c\delta_{x=0}\hbox{ in the sense of distributions on }\R^N\,.
\end{equation}
The constant $c$ in the r.h.s. of (\ref{Rsd}) is called {\it the residue of $f$ at $0$}
and denoted $\hbox{Res}_0f$. Equivalently, the residue of $f$ at $0$ can be defined by
the relation
\begin{equation}
\label{RsdInt}
\la f,\Phi\ra=\hbox{Res}_0f\,\tfrac1{|S^{N-1}|}\int_{\R^N}\Phi(x)\frac{dx}{|x|^N}\,,
\end{equation}
whenever $\Phi(x)=\phi(|x|)$ with $\phi\in C^\infty_c((0,+\infty))$.

Any $f\in\DD'(\R^N\setminus 0)$ which is homogeneous of degree $-N$ can be extended as a
homogeneous distribution $\dot{f}$ of degree $-N$ on $\R^N$ if and only if $\hbox{Res}_0f
=0$. For each $\chi\in C^\infty_c(\R^N)$, set $X(x)=\int_0^1\grad_x\chi(tx)dt$; one has
$\chi(x)=\chi(0)+x\cdot X(x)$ and $X\in C^\infty(\R^N)$. Given $f\in\DD'(\R^N\setminus 0)$
that is homogeneous of degree $-N$ with $\hbox{Res}_0f=0$ and $\phi\equiv\phi(|x|)$ in
$C^\infty_c(\R^N)$ such that $\phi\equiv 1$ near $0$, the linear functional
$$
\chi\mapsto \la f,(1-\phi)\chi\ra +\la (xf)^\cdot,\phi X\ra
$$
is a homogeneous extension of $f$ to $\R^N$. Two homogeneous extensions of $f$ may differ
by a multiple of $\delta_{x=0}$.

\end{document}